\documentclass[12pt]{amsart}
\begin{document}
\numberwithin{equation}{section}

\def\1#1{\overline{#1}}
\def\2#1{\widetilde{#1}}
\def\3#1{\widehat{#1}}
\def\4#1{\mathbb{#1}}
\def\5#1{\frak{#1}}
\def\6#1{{\mathcal{#1}}}

\def\C{{\4C}}
\def\R{{\4R}}
\def\N{{\4N}}
\def\Z{{\4Z}}

%tentative title
\title[Boundary jets]{Boundary jets of holomorphic maps between strongly pseudoconvex domains}
\author[F. Bracci, D. Zaitsev]{Filippo Bracci \and Dmitri Zaitsev}
%\thanks{}
\address{F. Bracci: Dipartimento di Matematica, Universit\`a di Roma  ``Tor Vergata'', Via della Ricerca Scientifica 1, 00133 Roma, Italy.}
\email{fbracci@mat.uniroma2.it}
\address{D. Zaitsev: School of Mathematics, Trinity College Dublin, Dublin 2, Ireland.}
\email{zaitsev@maths.tcd.ie} 

\maketitle
%\tableofcontents

%\def\Label#1{\label{#1}{\bf (#1)}~}
\def\Label#1{\label{#1}}

% Standard sets

\def\cn{{\C^n}}
\def\cnn{{\C^{n'}}}
\def\ocn{\2{\C^n}}
\def\ocnn{\2{\C^{n'}}}
\def\je{{\6J}}
\def\jep{{\6J}_{p,p'}}

% Abbreviations

\def\dist{{\rm dist}}
\def\const{{\rm const}}
\def\rk{{\rm rank\,}}
\def\id{{\sf id}}
\def\aut{{\sf aut}}
\def\Aut{{\sf Aut}}
\def\CR{{\rm CR}}
\def\GL{{\sf GL}}
\def\Re{{\sf Re}\,}
\def\Im{{\sf Im}\,}
\def\U{{\sf U}}

\def\codim{{\rm codim}}
\def\crd{\dim_{{\rm CR}}}
\def\crc{{\rm codim_{CR}}}

\def\phi{\varphi}
\def\eps{\varepsilon}
\def\d{\partial}
\def\a{\alpha}
\def\b{\beta}
\def\g{\gamma}
\def\G{\Gamma}
\def\D{\Delta}
\def\Om{\Omega}
\def\k{\kappa}
\def\l{\lambda}
\def\L{\Lambda}
\def\z{{\bar z}}
\def\w{{\bar w}}
\def\Z{{\1Z}}
\def\t{\tau}
\def\th{\theta}

\def\la{\langle}
\def\ra{\rangle}

\emergencystretch15pt \frenchspacing

\newtheorem{Thm}{Theorem}[section]
\newtheorem{Cor}[Thm]{Corollary}
\newtheorem{Pro}[Thm]{Proposition}
\newtheorem{Lem}[Thm]{Lemma}

\theoremstyle{definition}\newtheorem{Def}[Thm]{Definition}

\theoremstyle{remark}
\newtheorem{Rem}[Thm]{Remark}
\newtheorem{Exa}[Thm]{Example}
\newtheorem{Exs}[Thm]{Examples}

\def\bl{\begin{Lem}}
\def\el{\end{Lem}}
\def\bp{\begin{Pro}}
\def\ep{\end{Pro}}
\def\bt{\begin{Thm}}
\def\et{\end{Thm}}
\def\bc{\begin{Cor}}
\def\ec{\end{Cor}}
\def\bd{\begin{Def}}
\def\ed{\end{Def}}
\def\br{\begin{Rem}}
\def\er{\end{Rem}}
\def\be{\begin{Exa}}
\def\ee{\end{Exa}}
\def\bpf{\begin{proof}}
\def\epf{\end{proof}}
\def\ben{\begin{enumerate}}
\def\een{\end{enumerate}}

\section{Introduction}

The goal of this paper is to initiate a study of holomorphic
mappings $F$ between two domains $D$ and $D'$ in $\C^{n+1}$, sending
$D$ to a subset $F(D)\subset D'$ whose shape approximates $D'$ as
much as possible. It is known since Poincar\'e \cite{Po} and
subsequent work by Tanaka \cite{Ta}, Chern-Moser \cite{CM} and
Fefferman \cite{F} that in general, there does not exist any
biholomorphic maps between two given bounded strongly pseudoconvex
domains in $\C^{n+1}$ with $n\ge 1$
(see \cite{BER} for a survey of further results in this direction).
On the other hand, there are clearly many biholomorphic maps $F$ from $D$ to open subsets of
$D'$. {\em Can one impose any condition on $F$ making it ``close''
to be biholomorphic between $D$ and $D'$ without losing the
existence of such maps in general?} If one does not keep the
condition $F(D)\subset D'$, the Chern-Moser theory gives an estimate
on the maximal possible contact order between the boundaries of $D$
and $F(D)$. Our goal here is to study this question {\em under the
assumption} $F(D)\subset D'$.

We treat the problem locally at the given points $p\in\d D$ and $p'\in \d D'$ and consider the set
$\je_{p,p'}^0(D,D')$ of all germs $F$ at $p$ of holomorphic maps from $D$ to $D'$
sending $p$ to $p'$ in the non-tangential sense.
That is, an element in $\je_{p,p'}^0(D,D')$ is represented
by a holomorphic map $F\colon U\cap D \to D'$ with $U$ being a neighborhood of $p$,
such that $F(Z_n)\to p'$ whenever $Z_n\to p$ non-tangentially in $U\cap D$
(i.e.\ the distance from $Z_n$ to $p$ does not exceed its distance to $\d D$
times a constant multiple).
We shall assume both $D$ and $D'$ to be strongly pseudoconvex
with smooth boundaries and choose local holomorphic coordinates
$(z,w)$ and $(z',w')$ near $p$ and  $p'$ respectively
where $p=p'=0$ and $D$, $D'$ are locally given by
\begin{equation}\Label{init-norm}
\Im w>\|z\|^2+O(3), \quad \Im w'>\|z'\|^2+O(3),
\end{equation}
where $\|z\|:=|z_1|^2+\cdots+|z_n|^2$.
In order to speak about a contact order between $F(D)$ and $D'$
we need to introduce a differential of $F$ at the boundary point $p$,
which we again understand in the non-tangential sense (see \S{2} for precise definition).

The first question is what the possible non-tangential differentials
that may occur in this way are. As a first preliminary result we
give a complete characterization in terms of the singular values.
Recall that every complex $n\times n$ matrix $C$ admits its singular
value decomposition $C=U_1 D U_2$, where $U_1,U_2\in \U(n)$ are
unitary and $D$ is diagonal with real nonnegative entries
$\mu_1\ge\cdots\ge \mu_n\ge 0$ (this can be shown, {\sl e.g.}, by
using the polar decomposition $C=UH$ with $U$ unitary and $H$
hermitian and by further diagonalizing $H$). The entries of $D$ are
uniquely determined by $C$. We have the following characterization:

\bp\Label{firstcond}
A linear map $L\colon \C^{n+1}\to \C^{n+1}$ is
the non-tangential differential of a germ $F\in \je_{p,p'}^0(D,D')$
if and only if, in the chosen coordinates, it is of the form
\begin{equation}\Label{L}
L=\left(\begin{matrix} C & A \\ 0 & \l
\end{matrix}\right),
\end{equation}
where $\lambda>0$ is a real number, $A \in \cn$ is a complex vector
and $C$ is a complex $n\times n$ matrix whose singular values
$\mu_1\geq \ldots \geq \mu_n\geq 0$ satisfy $\mu_j \leq
\sqrt{\lambda}$ for $j=1,\ldots, n$.  \ep

In particular from \eqref{L} it follows  that a germ $F$ of
holomorphic map from $D$ to $D'$ sending $0$ to $0$ which is
non-tangentially differentiable at $0$ is a {\sl contact map}, in
the sense that its non-tangential differential  maps the complex
tangent space $T_0^c\d D$ into $T_0^c\d D'$. In case $D, D'$ are
bounded strongly pseudoconvex domains and $F:D\to D'$ is holomorphic
(i.e. it is not just a germ near a boundary point), this latter fact
follows also from Abate's generalization of the classical
Julia-Wolff-Carath\'eodory theorem (\cite{Ab}). Thus, in a certain
sense, Proposition~\ref{firstcond} can be interpreted as a
Julia-Wolff-Carath\'eodory theorem in the local.

If $F$ is as in Proposition \ref{firstcond}, set
$\alpha_j:=\mu_j/\sqrt{\lambda}$ for $j=1,\ldots, n$. We call the
numbers $1\ge \alpha_n \ge\ldots\ge \alpha_1\ge 0$ the {\em
singular values of the (non-tangential) differential of $F$ at $p$}.
It turns out that these numbers do not depend on the choice of coordinates
$(z,w)$ and $(z',w')$ provided \eqref{init-norm} holds
(see Lemma~\ref{firstcondbis}).
On the other hand, one can easily eliminate $A$
by composing $F$ with a suitable automorphism of the corresponding Siegel domain
$\Im w > \|z\|^2$ of the form
\begin{equation}\Label{auto}
g_a(z,w):=\frac{(z+aw, w)}{1 -2 i \la z,a \ra-i \|a\|^2 w}.
\end{equation}
Hence these are the only ``first order'' invariants of $F$ at $p$ and,
roughly speaking, they read the ratios of ``squeezing'' by $F$ in complex
tangent directions comparing to the normal direction.
The nearer to $1$ the singular values are,
the similar to $D'$ the image $F(D)$ looks like near $p'$.

As the next step, we study the conditions on the ``higher order jets'' of $F$ at $p$.
In order to make it meaningful to talk about jets at boundary points,
we shall assume $F$ to have smooth extension through the boundary.
That is, we consider the subset $\je_{p,p'}(D,D')$ of $\je_{p,p'}^0(D,D')$
consisting of all germs $F$ having representatives extending smoothly to
some neighborhoods of $p$.
It is not hard to see that if $\alpha_j<1$ for all $j$,
then there are no  restrictions on the set of possible higher order
jets of maps in $\je_{p,p'}(D,D')$ whose differentials at $p$ have
the given singular values $1>\alpha_1\ge\ldots\ge \alpha_n\ge 0$.
However, our Proposition~\ref{firstcond} above implies
that even the choice of $F$ with $1=\alpha_1 =\ldots = \alpha_n$
is always possible giving a better contact of $F(D)$ with $D'$.
Our next question is now to examine the possible restrictions
on the higher order jets of $F$ in this ``extreme case''.

Let $F\in \jep(D,D')$ be such that all the singular values of its
differential at $p$ are $1$. Then in view of the remarks above,
one can choose the coordinates $(z,w)$ and $(z',w')$ preserving \eqref{init-norm}
such that $dF_p$ becomes the identity $\id$.
The property of having $dF_p=\id$ in suitable coordinates where \eqref{init-norm} holds,
admits a natural higher order generalization:

\bd\Label{map-extrem}
A germ $F\in \jep(D,D')$ is said to be {\em $k$-flat}
if there exist local coordinates $(z,w)$ and $(z',w')$
vanishing at $p$ and $p'$ respectively, where the hypersurfaces
$\d D$ and $\d D'$ are respectively in their Chern-Moser normal forms
and such that $F=\id+o(k)$.
\ed

In other words, a germ $F\in \jep(D,D')$ is $k$-flat if and only if
there exist two local biholomorphic maps of $\C^{n+1}$, $h$ at $p$ and
$g$ at $p'$, such that $h(p)=g(p')=0$, $h(\partial D)$ and
$g(\partial D')$ are in their Chern-Moser normal forms and
$g\circ F\circ h^{-1}=\id+o(k)$, where  $H=o(k)$ means that $H$ and all its
derivatives of order less than $k+1$ are $0$ at $0$.
Thus a germ $F\in \jep(D,D')$ is {\em $1$-flat} if and only if all
singular values of $dF_p$ are $1$ and, as a consequence of
Proposition \ref{firstcond}, there always exist $1$-flat germ
for any $D$ and $D'$.

\br
It follows from the construction of the normal form in \cite{CM}
that the Chern-Moser normalizations in Definition~\ref{map-extrem}
are only needed to be chosen for the terms of weight $\le k$,
where as usual, the weight of $z$ and $\bar z$ is $1$ and of $w$ is $2$.
\er

Using Chern-Moser normal forms, we give
a complete description of the  second order jets for maps in
$\jep(D,D')$ whose first jet is the identity (Theorem \ref{secondcond}).
In particular it turns out that the space of possible second order jets
has its interior described by simple algebraic inequalities.
That is, for any $2$-jet in the interior, there exists
a germ $F\in\jep(D,D')$ with that jet and no further restrictions arise on the possible
jets of order three or higher.
This gives a more precise description of possible $1$-flat germs.

In contrast to $1$-flat germs, $2$-flat germs may not exist at all
for some $D$ and $D'$. This latter fact is somewhat related to the
rigidity phenomena for {\em self-maps} known as ``Burns-Krantz type
theorems'' (see \cite{B-K}, \cite{Hu}, \cite{BZZ}). We show that the
existence of $2$-flat germs implies a nontrivial geometric condition
on $D$ and $D'$ expressed as follows. We say that two real
hypersurfaces $M$ and $M'$ in $\C^{n+1}$ passing through a point $q$
are {\em tangent at $q$ up to weighted order $k$} if, for some (and
hence any) local defining function $\rho$ of $M'$ and some (and
hence any) local parametrization $\gamma\colon \cn\times\R\to
\C^{n+1}$ of $M$ with $\gamma(0)=q$ and $d\gamma_0(\cn\times\{0\})$
being the complex tangent space of $M$ at $q$, the composition
$\rho\circ \gamma$ vanishes at $0$ up to weighted order $k$, where
as before we assign weight $1$ to the coordinates in $z\in\cn$ and
weight $2$ to the coordinate in $u\in\R$. We now call $(\d D,p)$ and
$(\d D',p')$ {\em equivalent up to weighted order $k$} if there
exists a local holomorphic diffeomorphism of $\C^{n+1}$ near $p$,
sending $p$ to $p'$ and $\d D$ to another real hypersurface, which
is tangent to $\d D'$ up to weighted order $k$ at $p'$.

Our result for $2$-flat germs can now be stated as follows:

\bt\Label{nolabel} Let $D, D'\subset \C^{n+1}$ be two domains with
smooth boundaries such that $p\in \d D$, $p'\in\d D'$ and $\d D, \d
D'$ are strongly pseudoconvex at $p$ and $p'$ respectively. Then
there exist $2$-flat maps in $\jep(D,D')$ if and only if $(\d D,p)$ and
$(\d D',p')$ are equivalent up to weighted order $5$. \et

The outline of the paper is the following. In the second section we
prove the ``only if'' part of Proposition~\ref{firstcond} and
discuss the first jets. In the third section we recall the Chern-Moser
theory as needed for our purposes, describe the possible second jets for $1$-flat
germs (Theorem \ref{secondcond}), finish the proof of Proposition
\ref{firstcond}, give some equivalent conditions for $2$-flatness
and prove Theorem \ref{nolabel}. Finally, there is an Appendix where
we collected some auxiliary results needed in the various proofs.

\section{First order Jets}\Label{2}
Let $F\colon D\to \C^m$ be holomorphic for some $m$.
We say that $F$ is {\em non-tangentially differentiable} at $p$ if there exists
a point $p'\in\C^m$ and a linear map
$dF_p\colon \C^{n+1}\to \C^m$ such that
\begin{equation}\Label{diff-conv}
F(Z_k)=p'+ dF_p(Z_k-p)+o(\|Z_k-p\|), \quad k\to\infty,
\end{equation}
holds for any sequence $(Z_k)$ of points in $D$ converging non-tangentially to
$p$. We call $dF_p$ the {\em non-tangential differential} of $F$ at
$p$.

We shall consider the case when $D, D'$ are domains in $\C^{n+1}$ with smooth boundaries
and strongly pseudoconvex points $p\in \d D$, $p'\in \d D'$.
In the sequel, we shall assume $p=p'=0$
and choose local holomorphic coordinates
$(z,w)\in \cn \times \C$ and $(z',w')\in \cn \times \C$ vanishing at the origin
% (considering
%$D$ and $D'$ as domains in two distinct manifolds)
such that $D$ and $D'$ are locally given by
%expressions of the form
\begin{equation}\Label{normal-1}
D=\{\Im w >\|z\|^2+O(|(z,\Re w)|^3)\}, \quad D'=\{\Im w' >\|z'\|^2+O(|(z',\Re w')|^3)\}.
\end{equation}

We will obtain the ``only if'' statement of Proposition
\ref{firstcond} as consequence of the following  lemma:

\bl\Label{firstcondbis} Let $D, D'$ be two domains in $\C^{n+1}$
of the form \eqref{normal-1} and
$F\in\6J_{0,0}^0(D,D')$ be non-tangentially differentiable at $0$ with
differential $dF_0$. Then $dF_0$ is given by the block matrix
\begin{equation}\Label{matrixuno}
dF_0=\left(\begin{matrix} C & A \\ 0 & \l
\end{matrix}\right),
\end{equation}
where $\lambda>0$ is a real number, $A \in \cn$ is a complex vector
and $C$ is a complex $(n\times n)$-matrix whose singular values
$\mu_1\geq \ldots \geq \mu_n\geq 0$ satisfy $\mu_j \leq
\sqrt{\lambda}$ for $j=1,\ldots, n$. The ratios
$\mu_j/\sqrt{\lambda}$ for $j=1,\ldots, n$ are invariant under
coordinates changes preserving \eqref{normal-1}. \el

For the proof we need the following elementary result,
whose proof is supplied for the reader's convenience.

\bl\Label{preli}
Let $D$ be a domain in $\C^{n+1}$ having $0$ as a smooth boundary point
and $F\colon D\to \C^m$ be holomorphic.
If $dF_0$ is the non-tangential
differential of $F$ at $0$, then
$dF_{Z_k}\to dF_0$ for any sequence $(Z_k)$ in $D$ converging
non-tangentially to $0$. \el

\bpf Without loss of generality, $F(0)=0$ and $dF_0 = 0$. Let $
\{Z_k\}$ be any sequence in $D$ converging non-tangentially to $0\in
\partial D$. It suffice to show that $\frac{\d F}{\d Z^l}(Z_k)\to 0$
for every $l=1,\ldots,n+1$, where we use the notation
$Z=(Z^1,\ldots,Z^{n+1})$. We give a proof for $l=1$, the other cases
being completely analogous. Since $\{Z_k\}$ converges
non-tangentially in $D$, there exists $\eps>0$ such that, for any
other sequence $\{\2Z_k\}$ with $\|\2Z_k-Z_k\|\le \eps \|Z_k\|$, one
has $\2Z_k\in D$ for all $k$ and $\{\2Z_k\}$ also converges to $0$
non-tangentially in $D$. By the Cauchy Integral Formula, we have
\begin{equation}\Label{est}
\begin{split}
\left|\frac{\d F}{\d Z^1} (Z_k)\right|&= \left|\frac{1}{2\pi i}
\int_{|\zeta-Z_k^1|=\eps \|Z_k\|
}\frac{F(\zeta,Z_k^2,\ldots,Z_k^{n+1})}{(\zeta-
Z_k^1)^{2}}d\zeta\right| \\&\leq \frac{1}{\eps
\|Z_k\|}\max_{\zeta}|F(\zeta,Z_k^2,\ldots,Z_k^{n+1})|.
\end{split}
\end{equation}
It remains to choose $\zeta$ with $|\zeta-Z_k^1|=\eps \|Z_k\|$
such that the maximum in \eqref{est} is attained for $\2Z_k := (\zeta,Z_k^2,\ldots,Z_k^{n+1})$
and use \eqref{diff-conv} with $Z_k$ replaced by $\2Z_k$.
 \epf

\bpf[Proof of Lemma \ref{firstcondbis}] We first observe that $dF_0$
must send the upper half-space $\Im w \ge 0$ into itself. Otherwise
there would exist a non-tangentially convergent sequence $Z_k =
v\eps_k$, where $v\in \{(z,w) : \Im w > 0\}$ is a vector with
$dF_0(v)$ not contained in $\Im w \ge 0$ and $\{\eps_k\}$ a sequence
of positive real numbers converging to $0$. The latter would be in
contradiction with \eqref{diff-conv} and \eqref{normal-1}. Hence
$dF_0$ sends the real hyperplane $\Im w = 0$ into itself and, since
it is complex-linear in view of Lemma~\ref{preli}, also the complex
hyperplane $w = 0$ into itself. Putting everything together, we
conclude that $dF_0$ is of the form \eqref{matrixuno} with some
matrices $C$ and $A$ and a real number $\l\ge 0$.

The second step consists of showing that $\l>0$. Suppose, on the
contrary, that $\lambda=0$. Since $D'$ is strongly pseudoconvex
at $0$, it is easy to construct a continuous plurisubharmonic (peak)
function $\phi$ defined in a neighborhood of $0$ in $\C^{n+1}$
such that $\phi(0)=0$, $d\phi_0=-d(\Im w)$ and $\phi(Z)<0$ for
$Z\in \1{D'}\setminus\{0\}$. Furthermore, it is easy to
extend $\phi$ to a continuous plurisubharmonic function
$\psi$ defined on the whole $D'$ by setting
\begin{equation}
\psi(Z):=
\begin{cases}
\max(\phi(Z),-\eps) & \text{ for } \|Z\|<\delta \cr -\eps & \text{otherwise},
\end{cases}
\end{equation}
where $\delta>0$ and $\eps>0$ are chosen such that $\phi(Z)<-\eps$
for $Z\in D'$ with $\|Z\|=\delta$. Note that $\psi$ coincides with
$\phi$ in a neighborhood of $0$ in $\1{D'}$. Then $\l\ne 0$
follows from the Hopf lemma applied to $\psi\circ F$
restricted to a disk in the complex line $\{(z,w): z=0\}$
that is contained in $D$ and tangent to the boundary $\d D$ at $0$.

The third step is to show that the ratios $\mu_j/\sqrt{\lambda}$ do
not depend on the coordinates chosen and that the inequalities
$\mu_j/\sqrt\l\le 1$ hold. Using the singular value decomposition
of $C$, we can compose $F$ with suitable unitary linear transformations
of $\C^{n+1}$ on the right and on the left,
such that both forms \eqref{normal-1} are preserved
and $C$ becomes diagonal with real entries $\mu_1 \geq \ldots \mu_n \geq 0$
equal to its singular values.
Furthermore, composing with a dilation
$(z,w)\mapsto (\l z, |\l|^2w)$, we may assume that $\l=1$.

Now consider any changes of coordinates $(z,w)\mapsto \phi_1(z,w)$
and $(z',w')\mapsto \phi_2(z',w')$ preserving \eqref{normal-1}. Then
the differentials $(d\phi_1)_0$ and $(d\phi_2)_0$ must be of the
form
\[
\left(\begin{matrix}U_j\l_j & \ast \\
0 & \l_j^2
\end{matrix}\right),
\quad j=1,2,
\]
where $\l_j$'s are real positive and $U_j$'s are unitary.
Furthermore, in order to keep the above normalization of $dF_0$,
we must have $\l_1=\l_2$.
Then in these new coordinates, we have
\begin{equation}\Label{}
dF_0 = \left(\begin{matrix}U_2 C U_1^{-1} & \ast \\
0 & 1
\end{matrix}\right)
=\left(\begin{matrix}\2C & \ast \\
0 & 1
\end{matrix}\right),
\end{equation}
and it follows that the singular values of $\2C$ coincide with those of $C$.
This shows that the ratios $\mu_j/\sqrt{\l}$
are invariants of $dF_0$.

To show that $\mu_j/\sqrt{\lambda}\leq 1$ for all $j=1,\ldots, n$,
or $\mu_j$ in our normalization, it suffices to show that $\|C\|\le 1$.
By contradiction, suppose that
$\|C \xi\|>1$ for a vector $\xi\in\C^n$ with $\|\xi\|=1$.
We change local holomorphic coordinates in $\C^{n+1}$ near $0$
such that $\partial D$ and $\partial D'$ are approximated by the ball
$\{\|Z-(0,i/2)\|<1\}$ up to order $3$ at $0$.
Such coordinate change can be chosen to be the identity up to order $2$,
so that the matrix of $dF_0$ does not change.
Then, in view of Lemma~\ref{preli}, we can choose discs
\begin{equation}\Label{horiz-discs}
f_k\colon \D\to D, \quad f_k(\zeta):=\Big(\zeta\big(1-\frac{\eps}2\big)
\frac{\xi}{k}, \frac{i}{k^2}\Big)\in \C^n\times \C,
\end{equation}
where $\D$ is the unit disc in $\C$,
with $\eps>0$ sufficiently small such that
\begin{equation}\Label{schwarz}
\left\|\frac{1-\eps}{1+\eps}dF^z_{Z_k}(\xi)\right\|>1,
\end{equation}
where $F^z\in \C^n$ denotes the
tangential component of $F$ and $Z_k:=f_k(0)$. By the attraction
property (Lemma~\ref{attract}), for  $\eta:= \frac{1-\eps}{1-\eps/2}$, we may assume that
the images $F(f_k(\eta\D))$ are contained in a sufficiently small neighborhood of $0$.
Choose $r_k>0$ such that the central projection $\pi_k$ from $(0,-r_k)$ onto
the hyperplane $\{w=i/k^2\}$ sends the ball $\{\|Z-(0,i/2)\|<1\}$ in $\C^{n+1}$
into the ball with center $(0,i/k^2)$ and radius $\frac{1+\eps/2}k$ in the hyperplane $\{w=i/k^2\}$.
Then there exists sufficiently small neighborhood $U$ of $0$
such that, for $k$ sufficiently large,
$\pi_k$ sends $U\cap D'$ into the ball with center
$(0,i/k^2)$ and radius $(1+\eps)/k$.
Together with \eqref{schwarz}, we
reach a contradiction with the Schwarz lemma for
$\pi_k\circ F\circ f_k$ restricted to $\eta\D$.
 \epf

\br We say that a holomorphic map $F$ from $D$ into $D'$
which extends smoothly to $p$ and maps $p$ to $q$ and {\em sends
$\partial D$ into $\partial D'$ up to order $k$ at $p$} if, for
some (and hence any) local defining function $\rho_2$ of $\partial
D'$, its pullback $\rho_2\circ F$ is $o(k)$ on $\partial D$ ({\sl
e.g.}, $F$ sends $\partial D$ into $\partial D'$ up to order $1$
if $dF_p(T_p\partial D)\subseteq T_q\partial D'$). From the
previous discussion it is clear that if $F$ extends smoothly to
$p$ and sends $\partial D$ into $\partial D'$ up to order $2$,
then the singular values of $dF_p$ at $p$ are all equal to $1$.

\er

Let now $F\in\je_{0,0}^0(D,D')$. In the proof of
Lemma~\ref{firstcondbis} we have seen that in case all singular
values of the differential of $F$ at $0$ are equal to $1$,
we can choose coordinates $(z,w)$ and $(z',w')$
such that \eqref{normal-1} holds and
\[
d F_0=\left(\begin{matrix}\id  & A \\
0 & 1
\end{matrix}\right),
\]
for some complex vector $A\in \cn$. Using automorphisms $g_a$ of the
Siegel domain $\{(z,w) \in \cn \times \C: \Im w>\|z\|^2\}$ given by
\eqref{auto} we can replace $F$ by $g_{-A} \circ F$ to make its
differential at the origin equal to the identity. Here and in the
sequel we set
\begin{equation}\Label{braket}
\la z,\zeta \ra:= z_1\bar\zeta_1 + \ldots + z_n \bar\zeta_n.
\end{equation}
Note that $g_a$ preserves the form \eqref{normal-1}.

According to \cite{CM}, we can find germs of
biholomorphisms $h_1$ and $h_2$ such that $h_1(0)=h_2(0)=0$ and
$d(h_1)_0=d(h_2)_0=\id$ and $h_1(D), h_2(D')$ are in their Chern-Moser
normal forms. If $dF_0=\id$, we have $h_2\circ F \circ h_1^{-1}=\id+o(1)$.
Therefore we have:

\bc\Label{damettere} A germ $F \in \je_{p,p'}^0(D,D')$ is  $1$-flat if and
only if all singular values of its differential at $p$ are equal to
$1$. \ec

\section{Second order jets for $1$-flat maps}
As a matter of notations, for $m \in \N$, we use the symbol $O(m)$
to represent any (smooth) function which vanishes at the origin
together its derivatives up to order less than $m$. The symbol
$o(m)$ for $m \in \N$ means that also the $m$-th derivative is zero
at the origin. Whenever we need to state explicitly that a function
depending on several (complex or real) variables vanishes at the
origin together with all its partial derivatives with respect to a
certain variable---say $u$---up to the order $m$, we write such a
function as $O(u^m)$. Also we freely mix and add these notations.
For instance the function $3u^4v^3+v^2u^5$ can be written as $O(7)$,
or as $O(v^2)$ or $O(v^2)+O(u^4)$ or even as $O(u^4)+O(5)$. The same
notation is used for the small Landau's symbol $o$.

In this section we assume $F\in\jep(D,D')$ to be {\em $1$-flat}.
Recall that our notation $\jep(D,D')$ was reserved for the holomorphic map germs
having smooth extensions to some neighborhoods of $p$.
Arguing as in the previous section we may assume $p=p'=0$, $dF_0=\id$
and $\partial D, \partial D'$ given (locally) by expressions of the
form $\Im w=\|z\|^2+O(3)$. In order to simplify the notation, we use
the symbol $\je(D,D')$ to denote the germs of holomorphic maps from
$D$ to $D'$ which are smooth at $0$ and such that $F(0)=0$.

As a matter of notation, if $f\colon \C^{n+1}=\cn\times \C \to
\C^{n+1}$ is expandable at the origin, with homogeneous expansion
$f(z,w)=\sum_{\nu}f_\nu(z,w)$, we are going to denote by
$f^z_j(z,w)$ the projection to $\C^n_z$ of the homogeneous
polynomial vector $f_j$ and by $f_j^w$ the projection to $\C_w$.
Moreover, for a homogeneous polynomial $P(z,w)$ of degree $j$, we
write $P(z,w)=\sum_{\nu=0}^j P_{\nu, j-\nu}(z,w)$, for
$P_{l,k}(z,w)=C_{lk}(z)w^k$, where $C_{lk}(z)$ is a homogeneous
polynomial of degree $l$ in the $z$'s.

To deal with jets of order two we need however to have better
expansions for the normal forms of the domains.

Following \cite{CM}, we assign weight $1$ to $z_j, \bar z_j$
(for $j=1,\ldots, n$) and weight $2$ to $u=\Re w$. A real
polynomial $P(z, \bar z,u)$ is of {\em weighted degree} $m$ if it
is a linear combination of monomials of type $z_{j_1}\cdots
z_{j_k}u^l$ with $k+2l=m$. With this notation, Chern-Moser
normal forms for $\partial D$ and $\d D'$ can be written as
\begin{equation}\Label{moser}
\begin{split} &\d D=\{\Im w=\|z\|^2+\sum_{\mu\geq 4} \varphi_\mu(z,\bar z,\Re
w)\}, \\ &\d D'=\{\Im w=\|z\|^2+\sum_{\mu\geq 4} \varphi'_\mu(z,\bar
z,\Re w)\},
\end{split}
\end{equation}
where $\varphi_\mu$ and $\varphi'_\mu$  are
real weighted homogeneous polynomials of weighted degree $\mu$
which are linear combinations of monomials, each of which is
divisible by  $z_{j_1}z_{j_2}\bar z_{k_1}\bar z_{k_2}$ for some
$j_1,j_2,k_1,k_2\in \{1,\ldots, n\}$. In particular,
$\varphi_4,\varphi'_4, \varphi_5,\varphi'_5$ have no dependence in $\Re w$.
Also, $\hbox{tr}(\varphi_{4})\equiv 0$ and
$\hbox{tr}^2(\varphi_5)\equiv 0$ and similarly for $\phi'_4,\phi'_5$, where, if $Q_{(j,k)}(z,\bar
z)$ is a  polynomial of degree $j$ in $z_\alpha$ and degree $k$ in
$\bar z_\beta$ given by
\[
Q_{(j,k)}(z,\bar z)=\sum a_{\alpha_{1}\ldots \alpha_{j}
\beta_{1}\ldots\beta_{k}}z_{\alpha_{1}}\cdots z_{\alpha_{j}}\bar
z_{\beta_{1}}\cdots \bar z_{\beta_{k}},
\]
with the $a_{\alpha_{1}\ldots \alpha_{j}
\beta_{1}\ldots\beta_{k}}$'s symmetric with respect to $\alpha_{1}
\ldots, \alpha_{j}$ and respect to $\beta_{1}, \ldots, \beta_{k}$,
then
\[
\hbox{tr} (Q_{(j,k)}):=\sum_{\alpha_{1},\ldots, \alpha_{j-1}
,\beta_{1},\ldots, \beta_{k-1}}\left(\sum_{\alpha_{j}=\beta_{k}}
a_{\alpha_{1}\ldots \alpha_{j} \beta_{1}\ldots\beta_{k}}\right)
z_{\alpha_{1}}\cdots z_{\alpha_{j-1}} \bar z_{\beta_{1}}\cdots
z_{\beta_{k-1}}.
\]
Actually a Chern-Moser normal form as defined in~\cite{CM}
involves further trace conditions on higher order terms
that we won't need here.
Notice that the Chern-Moser normal form of a domain is not unique,
but it is parametrized by the automorphisms  of the quadric
$\{\Im w=\|z\|^2\}$ fixing the origin.

\bt \Label{secondcond}
Let $D$ and $D'$ be in their Chern-Moser normal
forms \eqref{moser} and $F\in \6J_{0,0} (D,D')$ with $d F_0=\id$. Then
\begin{equation}\Label{secf}
F(z,w)=(z,w)+
\big(F^z_{1,1}(z,w)+F^z_{0,2}(w),F^w_{0,2}(w)\big)+O(3)
\end{equation}
with
\begin{equation}\Label{ord}
\begin{split}
&\Im F^w_{0,2}(1)\geq 0,\\ &\left[\Re \la z,F^z_{1,1}(z,1)\ra -
\|z\|^2 \Re F^w_{0,2}(1)\right]^2\\ &\leq \Im F^w_{0,2}(1)
\left[\varphi_{4}(z)-\varphi'_{4}(z)-2
\|z\|^{2}\Im \la z, F_{1,1}^{z}(z,1)\ra-\|z\|^{4}\Im
F^w_{0,2}(1)\right]
\end{split}
\end{equation}
and one has $F_{3,0}^w \equiv 0$ and
\begin{equation}\Label{ord2}
\varphi_{4}(z)-\varphi'_{4}(z)-2 \|z\|^{2}\Im \la z,
F_{1,1}^{z}(z,1)\ra\geq 0.
\end{equation}
 On the other hand, for any
choice of the 2nd order terms in \eqref{secf}
satisfying~\eqref{ord} with  strict inequalities for all $z\in
\C^{n}\setminus\{0\}$ there exists  $F\in \6J_{0,0}(D,D')$ of
the form \eqref{secf}. \et

\bpf We shall use Corollary~\ref{weights} applied to the basic
condition $F(D)\subseteq D'$. In view of~\eqref{moser}, a
parametrization for $\partial D$ is given by
\begin{equation}\Label{parametro} \cn
\times \R \ni (z,u) \mapsto Z=\big(z, u+i\|z\|^2+i\sum_{\mu\geq
4}\varphi_\mu^1(z,\bar z,u)\big).
\end{equation}
Therefore the basic condition
becomes
\begin{equation}\Label{basic}
\begin{split}
\sum_{\mu\geq 4}&\varphi_\mu(z,\bar{z},u)+\sum_{k\geq 2}\Im F_{k}^{w}(Z)\\
&\geq 2\sum_{k\geq 2}\Re \la z, F_{k}^{z}(Z) \ra
+\|\sum_{k\ge2}F_{k}^{z}(Z)\|^{2}+\sum_{\mu\geq 4}\varphi'_\mu
(F^z(Z),\overline{F^z(Z)},\Re F^w(Z)),
\end{split}
\end{equation}
where $Z$ is as in \eqref{parametro} and $F_k$ denotes the component of weight $k$. Expanding
\eqref{basic} up to weighted order two and applying
Corollary~\ref{weights} we have
\begin{equation}\Label{order2}
\Im F_{2,0}^w\geq 0.
\end{equation}
Since $z\mapsto F_{2,0}^{w}(z)$ is holomorphic this means that
$F^{w}_{2,0}\equiv 0$.

Now expanding \eqref{basic} up to weighted order three and applying Corollary~\ref{weights} yields
\begin{equation}\Label{order3}
   \Im F^{w}_{3,0}(z)+\Im F^{w}_{1,1}(z,u+i\|z\|^{2})\geq 2 \Re\la z,
   F^{z}_{2,0}(z)\ra.
\end{equation}
Separating into terms of different bi-degree types and again using
Corollary~\ref{weights} we obtain two inequalities, namely
\begin{equation}\Label{p1}
\Im F^{w}_{1,1}(z,u)\geq 0
\end{equation}
and
\begin{equation}\Label{p2}
    \Im F^{w}_{3,0}(z)+\Im F^{w}_{1,1}(z,i\|z\|^{2})\geq 2\Re\la z,F_{2,0}^{z}(z)\ra.
\end{equation}
Inequality~\eqref{p1} is indeed an equality because
$F^{w}_{1,1}(z,u)$ is linear in $u$ and, since $z\mapsto
F^{w}_{1,1}(z,u)$ is holomorphic for any fixed $u$, it follows
furthermore that $F^{w}_{1,1}\equiv 0$. Now applying
Lemma~\ref{types} to \eqref{p2} we obtain $\Im F^{w}_{3,0}\equiv
0$ and hence $F^{w}_{3,0}\equiv 0$ for $z\mapsto
F^{w}_{3,0}(z)$ is holomorphic, and, consequently, $\Re\la
z,F_{2,0}^{z}(z)\ra\equiv 0$ for all $z\in \C^{n}$. This last
equality clearly implies $F^{z}_{2,0}\equiv 0$.

Therefore $1$-flatness implies that all terms of weighted order two
and three in the expansion of \eqref{basic} are zero.
Now we pass to  the weighted order four:
\begin{multline}\Label{order4}
    \varphi_{4}(z,\bar z)+\Im F^{w}_{4,0}(z)
+\Im F^{w}_{2,1}(z,u+i\|z\|^{2})+\Im F^{w}_{0,2}(u+i\|z\|^{2})\\ \geq 2\Re \la z, F^{z}_{3,0}(z)\ra
+ 2 \Re\la z, F^{z}_{1,1}(z,u+i\|z\|^{2} )\ra +\varphi'_{4}(z,\bar z).
\end{multline}
By Corollary~\ref{weights}, looking at terms of the lowest degree
in $(z,\bar z)$ we obtain that $\Im F^{w}_{0,2}(u)\geq 0$ and,
since the dependence on $u$ is quadratic, this is equivalent to
$\Im F^{w}_{0,2}(1)\geq 0$.

Now we can set $u=t\|z\|^2$ with $t\in \R$ in~\eqref{order4}, and
apply Remark~\ref{types2} to  terms of bi-degree $(2,2)$ in
$(z,\bar z)$:
\begin{multline}\Label{ww1}
t^2\|z\|^4 \Im F^w_{0,2}(1) +2t( \|z\|^4 \Re F^w_{0,2}(1)- \|z\|^2 \Re\la z, F^z_{1,1}(z,1)\ra)\\
+ \varphi_4(z,\bar z)-\varphi'_4(z,\bar z) -2\|z\|^2\Im \la z, F^z_{1,1}(z,1)\ra
-\|z\|^4 \Im F^w_{0,2}(1)\geq 0.
\end{multline}
For $z\neq 0$ fixed,  the left-hand side of \eqref{ww1} must be
greater than or equal $0$ for all $t$,
which is equivalent to \eqref{ord} and \eqref{ord2}
(for $\Im F^w_{0,2}(1)\ne0$, \eqref{ord2} follows from \eqref{ord}).

Finally, if both inequalities in~\eqref{ord} are strict for any
$z\neq 0$, then the lowest weighted order  nontrivial homogeneous
term in \eqref{basic} is positive for $(z,u)\ne0$
if we choose $F$ to be of the form \eqref{secf} without
higher order terms. Therefore \eqref{basic} will
always hold in a neighborhood of the origin.
This proves the last statement. \epf

\br\Label{less-1} It is apparent from the proof that the conclusions
of Theorem~\ref{secondcond} still hold  with $D, D'$ being given
in Chern-Moser normal forms only {\sl up to weighted order $5$} at $0$.
%{\sl i.e.} of the form  $\Im w>\|z\|^2+\varphi_4+O(\|z\|^5+|\Re w|\|z\|^3)\}$
%with $\hbox{tr}(\varphi_4)\equiv 0$.
\er

Now we are in the position to end the proof of Proposition
\ref{firstcond}:

\begin{proof}[End of the Proof of Proposition \ref{firstcond}] In
order to complete the proof of Proposition \ref{firstcond} we need
to show that given a matrix $L$ as in \eqref{L},
there exists $F\in \je_{0,0}(D,D')$ such that
$dF_0=L$. Using transformations tangent to $\id$ we can suppose that
$D, D'$ are in their Chern-Moser normal form (at least up to weighted
order four). Finally, acting with an automorphism \eqref{auto}
and a dilation $(z,w)\mapsto (\l z,|\l|^2w)$ on the left and with
unitary transformations $(z,w)\mapsto (Uz,w)$ on both sides,
we can reduce the general case to that of
\[
L=\left(\begin{matrix} \Delta & 0\\ 0 & 1
\end{matrix}\right),
\]
where $\Delta$ is a diagonal matrix with entries the singular values
$\alpha_j$'s, with $0\leq \alpha_n\leq\ldots\leq \alpha_1\leq 1$.
Now the argument is similar to that in the proof of Theorem
\ref{secondcond}. We look for $F$ of the form
$F(z,w)=(\Delta z, w)+\sum_{k\geq 2}(F^z_k(z,w), F^w_k(z,w))$ and
impose the condition $F(\d D)\subseteq \d D'$. Parametrizing $\d D$
with \eqref{parametro} we obtain
\begin{equation*}
\begin{split}
\|z\|^2&+\sum_{\mu\geq 4}\phi_\mu(z,\bar z, w)+\sum_{k\geq 2}\Im
F_k^w(Z)\\&\geq \|\Delta z\|^2+2\sum_{k\geq 2}\Re \la \Delta z,
F_{k}^{z}(Z) \ra +\|\sum_{k\ge2}F_{k}^{z}(Z)\|^{2}+\sum_{\mu\geq
4}\varphi'_\mu (F^z(Z),\overline{F^z(Z)},\Re F^w(Z)).
\end{split}
\end{equation*}
If all entries in $\Delta$ are $<1$, we choose $F_{0,2}^w$ with
$\Im F_{0,2}^w(1)>0$ and $F_k^z=F_{k+1}^w=0$ for $k\geq 2$, then it
follows that $F\in \je (D,D')$, that $F(\d D)\subset \d D'$ near $0$
and that $F$  has the required differential at $0$.

Let $l\leq n$ and suppose that $\alpha_1=\ldots=\alpha_l=1$ and
$\alpha_{k}<1$ for $k>l$. Let us write $z=(z',z'')\in \C^l\times
\C^{n-l}$. Also, with obvious meaning, write
$F^z_k=(F^{z'}_k,F^{z''}_k)$. Set $F^{z''}_k=0$ for $k\geq 0$.
The last statement in Theorem \ref{secondcond} gives sufficient conditions for the map
$(z',w)\mapsto F(z',0,w)$ to send the domain $D\cap \{z''=0\}$ into
the domain $D'\cap \{z''=0\}$.
Then the appropriate choice of $F^{z'}_k$ together with $\Im F_{0,2}^w(1)>0$
and the inequality $\|(0,z'')\|^2> \|\Delta (0,z'')\|^2$ for $z''\ne 0$
guarantees that $F$ sends $\d D$ into $\bar D'$ near $0$.
\end{proof}

\bl\Label{2-extremal-lem} Let $D,D'\subset\C^n$ be in Chern-Moser
normal forms and $F\in \6J_{0,0} (D,D')$ be with $dF_0=\id$ and
$\Im \la z, F^z_{1,1}(z,1) \ra \equiv0$. Then
\begin{equation}\Label{3-13}
\varphi_4(z,\bar z)\equiv \varphi'_4(z,\bar z),
\quad \varphi_5(z,\bar z)\equiv \varphi'_5(z,\bar z)
\end{equation}
and
\begin{equation}\Label{3-f1}
F(z,w)=(z,w)+(F^z_{1,1}(z,w),F^w_{0,2}(w))+(F^z_{1,2}(z,w)+F^z_{0,3}(w),F^w_{0,3}(w))+O(4),
\end{equation}
where
\begin{equation}\Label{nowork1}
\Im F^w_{0,2}(1)=\Im F^w_{0,3}(1)=0, \quad \|z\|^2\Re F^w_{0,2}(1)\equiv \Re \la z,
F^z_{1,1}(z,1)\ra
\end{equation}
and, for any $z$,
\begin{equation}\Label{nowork}
BC\ge A^2, \quad B\ge0, \quad C\ge0,
\end{equation}
where
\begin{equation}\Label{AB}
\begin{split}
A:=& -4\|z\|^4\Im \la z, F^z_{1,2}(z,1)\ra
+ \|z\|^2\phi_{2,2,1}(z,\bar z,1) - \|z\|^2\phi'_{2,2,1}(z,\bar z,1),\\
B:=& 3\|z\|^6\Re F^w_{0,3}(1) -2\|z\|^4\Re\la z, F^z_{1,2}(z,1)\ra
- \|z\|^4 \|F^z_{1,1}(z,1)\|^2, \\
C:=& -\|z\|^6\Re F^w_{0,3}(1) +2\|z\|^4\Re\la z, F^z_{1,2}(z,1)\ra
- \|z\|^4 \|F^z_{1,1}(z,1)\|^2\\
&+ \phi_{3,3,0}(z,\bar z) - \phi'_{3,3,0}(z,\bar z).
\end{split}
\end{equation}
Moreover $F^z_{4,0}\equiv0, F^w_{5,0}\equiv0, F^w_{3,1}\equiv0$.
 \el

\begin{proof}
We first prove that $\Im F^{w}_{0,2}(1)= 0$.
Indeed, we have $\Im F^{w}_{0,2}(1)\ge 0$ by Theorem~\ref{secondcond}.
If we had $\Im F^w_{0,2}(1)>0$, then dividing both sides of
\eqref{ord} by $\Im F^w_{0,2}(1)$ we would obtain that
\begin{equation}\Label{batistuta}
\varphi_{4}^{1}(z)-\varphi_{4}^{2}(z)-2 \|z\|^{2}\Im \la z,
F_{1,1}^{z}(z,1)\ra-\|z\|^{4}\Im F^w_{0,2}(1)\geq 0,
\end{equation}
and therefore $\varphi_{4}^{1}(z)-\varphi_{4}^{2}(z)\geq 0$
(since $\Im \la z, F^z_{1,1}(z,1) \ra \equiv0$).
Since $\hbox{tr} (\varphi_{4}^{1}(z)-\varphi_{4}^{2}(z)) =0$,
 the function $\varphi_{4}^{1}(z)-\varphi_{4}^{2}(z)$ is harmonic
and hence would be identically zero by the maximum principle.
Then \eqref{batistuta}
would imply $\Im F^w_{0,2}(1)= 0$, a contradiction.
Therefore $\Im F^w_{0,2}(1)=0$ and also $\varphi_{4}^{1}(z)-\varphi_{4}^{2}(z)\equiv 0$.
Hence \eqref{ord} implies  $\|z\|^{2}\Re F^{w}_{0,2}(1)\equiv
\Re \la z, F_{1,1}^{z}(z,1)\ra$.
This proves \eqref{nowork1}, except for $F^{w}_{0,3}$.

Summarizing, we have the equality in \eqref{ww1} for all $z$ and $t$.
Hence, by Lemma \ref{types}, we must have the equality also in \eqref{order4}.
In particular, separating types, we obtain the vanishing of
\begin{equation}\Label{j1}
 F^{w}_{4,0}, \quad F^{w}_{2,1}(z,1), \quad  F^{z}_{3,0}(z),
 \quad \|z\|^{2}\Re F_{2,1}^{w}(z,1)- 2\Re \la z, F^{z}_{3,0}(z)\ra.
\end{equation}

Hence in \eqref{basic} all terms of weighted order less or equal
to $4$ cancel each other, and we obtain the following inequality
for the terms of weighted order $5$:
\begin{multline}\Label{order5}
   \varphi_{5}(z, \bar z)+\Im F^{w}_{5,0}(z)+\Im F^{w}_{3,1}(z,w)+\Im F^{w}_{1,2}(z,w)\\
   \geq 2 \Re \la z, F_{4,0}^{z}(z)\ra + 2\Re \la z,
   F^{z}_{2,1}(z,w)\ra+2\Re \la z, F^{z}_{0,2}(w)\ra +\varphi'_{5}(z, \bar z),
\end{multline}
where the terms are evaluated at $(z,w)=(z,u+i\|z\|^2)$.
By Corollary~\ref{weights}, we can pass to the reduced inequality
involving only terms of degree $0$ in $u$. These are homogeneous
polynomials of the odd degree $5$ in $(z, \bar z)$ and
hence we have the  equality:
\begin{multline}\Label{u-zero}
\varphi_5(z, \bar z)+\Im F^w_{5,0}(z)+\Im F^w_{3,1}(z,
i\|z\|^2)+\Im F^w_{1,2}(z,i\|z\|^2) \\ = 2\Re \la z,
F^z_{4,0}(z)\ra + 2\Re \la z, F_{2,1}^z(z,i\|z\|^2)\ra +2\Re \la
z, F^z_{0,2}(i\|z\|^2)\ra +\varphi'_5(z, \bar z).
\end{multline}
Separating types we obtain
\begin{equation}\Label{5-0}
\Im F^w_{5,0}(z)=0,
\end{equation}
\begin{equation}\Label{4-1}
\Im F^w_{3,1}(z,i\|z\|^2)=2\Re\la z, F^z_{4,0}(z)\ra,
\end{equation}
\begin{equation}\Label{3-2}
\varphi_5(z, \bar z)+\Im F^w_{1,2}(z,i\|z\|^2)=2\Re\la z,
F^z_{0,2}(i\|z\|^2)\ra+2\Re\la z, F^z_{2,1}(z,i\|z\|^2)\ra +
\varphi'_5(z, \bar z).
\end{equation}
 Similarly we can repeat the argument for the terms of degree
$1$ in $u$ in~\eqref{order5} and separating types we obtain
\begin{equation}\Label{3-0}
\Im F^w_{3,1}(z,u)=0,
\end{equation}
\begin{equation}\Label{2-1}
u\|z\|^2 \Re F^w_{1,2}(z,1)=u \Re \la z, F^z_{2,1}(z,1)\ra +
2\|z\|^2u \Im \la z, F^z_{0,2}(1)\ra.
\end{equation}
Finally, repeating the process for the terms of degree $2$ in $u$
in~\eqref{order5}, we obtain
\begin{equation}\Label{u-2}
\Im F^w_{1,2}(z,u)=2\Re \la z, F^z_{0,2}(u)\ra.
\end{equation}
Now \eqref{5-0} and \eqref{3-0} imply  that $F^w_{5,0}\equiv 0$
and $F^w_{3,1}\equiv 0$. Then \eqref{4-1} gives $F^z_{4,0}\equiv
0$. From~\eqref{u-2}, dividing both sides by $u^2$ and noticing that
both maps $z\mapsto F^w_{1,2}(z,1)$ and $z\mapsto \la z,
F^z_{0,2}(1)\ra$ are holomorphic, we obtain
\begin{equation}\Label{f-1}
F^w_{1,2}(z,1)\equiv 2i\la z, F^z_{0,2}(1)\ra.
\end{equation}
Now substituting \eqref{f-1} into \eqref{3-2} (taking into account that
$F^w_{1,2}(z,i\|z\|^2)=-\|z\|^4 F^w_{1,2}(z,1)$ and
$F^z_{0,2}(i\|z\|^2)=-\|z\|^4 F^z_{0,2}(1)$), we obtain
\begin{equation}\Label{ci-1}
\varphi_5(z,\bar z)=\varphi'_5(z,\bar z)+2\|z\|^2\Im \la z,
F^z_{2,1}(z,1)\ra.
\end{equation}
Also, from \eqref{2-1} and \eqref{f-1} we find
\begin{equation}\Label{ei-1}
\Re \la z, F^z_{2,1}(z,1)\ra = -4 \|z\|^2 \Im \la z, F^z_{0,2}(1)\ra.
\end{equation}
Separating types, this means
$$\la F^z_{2,1}(z,1),z\ra = 4i \|z\|^2  \la z, F^z_{0,2}(1)\ra,$$
which, together with \eqref{ci-1} implies
\begin{equation}\Label{t-5}
\varphi_5-\varphi'_5=-8\|z\|^4 \Re \la z, F^z_{0,2}(1)\ra.
\end{equation}
Note that we have $\hbox{tr}^2(\varphi_5)\equiv\hbox{tr}^2(\varphi'_5) \equiv 0$.
Therefore using the uniqueness of the trace decomposition
(see \cite{CM}) we
conclude that $\varphi_5\equiv\varphi'_5$ and $F^z_{0,2}(1)=0$,
and hence $F^w_{1,2}(z,1)\equiv 0$ in view of \eqref{f-1}

It remains to show that $\Im F_{0,3}^w(1)=0$ and the inequalities in \eqref{nowork}.
We now can pass to the weighted order $6$ inequality, which yields
\begin{multline}\Label{order6}
\varphi_{6}(z, \bar z,u)+\Im F^{w}_{6,0}(z)+\Im
F^{w}_{4,1}(z,w)+
\Im F^{w}_{2,2}(z,w) + \Im F^{w}_{0,3}(w) \\
   \geq  2 \Re\la z, F_{5,0}^{z}(z)\ra +
2\Re\la z, F^{z}_{3,1}(z,w) \ra+ 2\Re\la z, F^{z}_{1,2}(z,w) \ra +
\|F^z_{1,1}(z,w)\|^2 + \varphi'_{6}(z, \bar z,u),
\end{multline}
where the terms are evaluated at $(z,w)=(z,u+i\|z\|^2)$.
Using Corollary~\ref{weights} for the terms of degree $0$ in $(z,\bar z)$,
we have $u^3\Im F^w_{0,3}(1)\geq 0$ which implies $\Im F^w_{0,3}(1)=0$.
Now we let $u=t \|z\|^2$ for
$t\in \R$ and look at the weighted order $6$ inequality using
Lemma~\ref{types} to pass to the terms of type $(3,3)$ in $(z,\bar z)$:
\begin{equation}\Label{reduced}
\begin{split}
\phi_{3,3,0}(z,\bar z) &+\phi_{2,2,1}(z,\bar z,1)t\|z\|^2 + \Im [F^w_{0,3}(1)(t+i)^3] \|z\|^6 \\
\geq &2\Re \la z, F^{z}_{1,2}(z,1) (t+i)^2\|z\|^4\ra
+\|F^z_{1,1}(z,1)(t+i)\|^2\|z\|^4\\+&\phi'_{3,3,0}(z,\bar z)
+\phi'_{2,2,1}(z,\bar z,1)t\|z\|^2 .
\end{split}
\end{equation}
In view of $\Im F^w_{0,3}=0$,  \eqref{reduced} leads to the
quadratic inequality $Bt^2 + 2At + C\ge0$ for all $t$ and $z$ with
$A,B,C$ as in \eqref{AB}. The latter inequality is clearly
equivalent to \eqref{nowork}.\end{proof}

\br\Label{weight-order} Observe that, for any $k$, the property that
one has the equality in \eqref{basic} up to weighted order $k$ does not
depend on the choice of coordinates. Indeed, \eqref{basic} is
obtained by substituting the parametrization
$$\gamma\colon(z,u)\mapsto F(z,u+i\|z\|^2+i\sum\phi_\mu(z,\bar z, u))$$
of $F(\d D)$ into the defining function
$$\rho(z,w):= \Im w - \|z\|^2 - \sum\phi'_\mu(z,\bar z, \Re w)$$
 of $\d D'$.
Then the equality in \eqref{basic} up to weighted order $k$
means that $\rho\circ \gamma$ vanishes up to weighted order $k$ at $0$.
Now we claim that
for any smooth defining function $\2\rho$ of $\d D'$
and any smooth parametrization $\2\gamma(\2z,\2u)=\gamma(z(\2z,\2u),u(\2z,\2u))$ of
$F(\d D)$ with $\frac{du}{d\2z}(0)=0$,
the weighted vanishing orders of $\2\rho\circ \2\gamma$ (in $(\2z,\2u)$)
coincides with that of $\rho\circ \gamma$  (in $(z,u)$).
Indeed, we have $\2\rho = \rho \alpha$ for a suitable function $\alpha$
and hence the weighted vanishing order of $\2\rho\circ \gamma$
is at least as high as that of $\rho\circ \gamma$.
Furthermore, writing $(z,u)=(A\2z+B\2u,C\2u) + O(\|\2z\|^2+ \2u^2)$
with suitable matrices $A,B,C$,
we see that also the weighted vanishing order of $\2\rho\circ \2\gamma$
is at least as high as that of $\rho\circ \gamma$.
Reversing the argument, we see that both vanishing orders are equal as claimed.
\er

We shall say that $F(\d D)$ is tangent to $\d D'$ at $0$ up to
weighted order $k$ if we have the equality in \eqref{basic} up to
weighted order $k$. The latter property is well-defined and does not
depend on coordinate choices in view of Remark~\ref{weight-order}.

\bp\Label{equivcond} Let $D,D'\subset\C^n$ be in their Chern-Moser
normal forms and $F\in \6J_{0,0} (D,D')$ be of the form
\eqref{secf}. The following conditions are equivalent:
\begin{enumerate}
\item The germ $F$ is $2$-flat (in the sense of Definition \ref{map-extrem});
\item $F(\d D)$ is tangent to $\d D'$ at $0$ up to weighted order $4$;
\item $\Im \la z, F^{z}_{1,1}(z,1) \ra\equiv0$.
\end{enumerate}
 \ep

\bpf Suppose that $F$ is $2$-flat and choose coordinates according
to Definition \ref{map-extrem} such that $F=\id + O(3)$. By
Theorem~\ref{secondcond}, $F_{3,0}^w\equiv0$ and therefore we have
the equality in \eqref{basic} up to weighted order $3$. Next, examining
terms of weighted order $4$ of types $(2,2,0)$ and $(0,0,2)$ in
$(z,\bar z,u)$ in \eqref{basic} we see that they only involve the
second derivatives of $F$ and $\phi_4^1-\phi_4^1$ (cf. \eqref{ww1}),
where the latter vanishes by Lemma~\ref{2-extremal-lem}
since $F_{1,1}^z\equiv 0$. Hence, by
Lemma~\ref{types}, the whole weighted homogeneous part of
\eqref{basic} of order $4$ must vanish. Thus (1) implies (2).

Now assume (2). In particular, we have the equality in \eqref{order4}
which, for the terms of type $(0,0,2)$ in $(z,\bar z,u)$ yields $\Im F_{0,2}^w(1)=0$.
Then the equality in \eqref{ww1} together with the trace decomposition implies (2) as in \cite{CM}.

Finally, assuming (3), applying Lemma~\ref{2-extremal-lem} and
arguing as before, we obtain (2) proving that (2) and (3) are in
fact equivalent. Now consider the parabolic automorphism of type
\[
g_r(z,w)=\frac{(z,w)}{1-rw}
\]
with $r=-\Re F^w_{0,2}(1)$. As shown in \cite{CM}, there exists a
unique transformation $h$ such that $h(0)=0$, $d_0h=\id$ and $\Re
h^w_{0,2}(1)=0$ and $\2D':=h(g_r(D'))$ is in its Chern-Moser normal
form. Then the map  $\tilde{F}=h \circ g_r \circ F$ satisfies (2)
(with respect to $D$ and $\2D'$) and, moreover, $\Re
\tilde{F}^w_{0,2}(1)=0$. As we have seen, (2) implies (3) and
therefore we can apply Lemma~\ref{2-extremal-lem} (identity \eqref{nowork1}) to $\2F$ to
conclude that $\2F=\id + O(3)$. Hence (1) holds as desired. \epf

\bpf[Proof of Theorem~\ref{nolabel}]
By definition, if there exists
a $2$-flat map $F\in \6J_{p,p'}(D,D')$, we have $F=\id+O(3)$ with respect
to some Chern-Moser normal coordinates for $\d D$ and $\d D'$
vanishing at $p$ and $p'$ respectively. Then by
Lemma~\ref{2-extremal-lem} (identity \eqref{3-13}),
the Chern-Moser normal forms of $\d D$ and
$\d D'$ coincide up to weighted order $5$ and therefore $(\d D,p)$ and $(\d D',p')$
are equivalent up to weighted order $5$.

Conversely, suppose $(\d D,p)$ and $(\d D',p')$ are
biholomorphically equivalent up to weighted order $5$.
Then it follows from the construction of the normal form in \cite{CM}
that there exist Chern-Moser normal forms for $\d D$ and $\d D'$ that
coincide up to weighted order $5$.
We will construct a map $F\in \6J(D,D')$ with $F=\id + O(3)$ of the form
\begin{equation}\Label{form}
F(z,w)= (z,w) + (\l_1 w^2 z, \l_2 w^3 + i\l_3 w^4) \quad
\end{equation}
with $\l_1,\l_2,\l_3$ being real numbers to be suitably chosen.
We first remark that with this choice of $F$ one always has the equality in \eqref{basic} up to
weighted order $5$. We now consider the corresponding inequality for
the terms of weighted order $6$:
\begin{equation}\Label{order6'}
\varphi^{1}_{6}(z, \bar z,u)+ \l_2 \Im (u+i\|z\|^2)^3 \\
\geq 2\l_1\|z\|^2 \Re(u+i\|z\|^2)^2 + \varphi_{6}^{2}(z, \bar z,u),
\end{equation}
which is equivalent to
\begin{equation}\Label{eq}
u^2\|z\|^2 (3\l_2-2\l_1) + \|z\|^6(-\l_2+2\l_1) \ge
\varphi_{6}^{2}(z, \bar z,u) - \varphi_{6}^{1}(z, \bar z,u),
\end{equation}
where
$\varphi_{6}^{2}(z, \bar z,u) - \varphi_{6}^{1}(z, \bar z,u)=O(\|z\|^6+u^2\|z\|^2)$.
Therefore we can choose $\l_1,\l_2$ to have the strict inequality in \eqref{eq}
whenever $z\ne0$.
We still have the equality for $z=0$, $u\ne 0$ and hence have to pass
to higher order terms to obtain strict inequality for all $(z,u)\ne 0$.
After further inspection of the terms of weighted order $7$ and $8$
we see that each of them, except $\l_3 u^4$,
is $o(\|z\|^6+u^2\|z\|^2)$ as $(z,u)\to 0$
due to the Chern-Moser normalization of the terms $\phi_\mu^j$.
Hence, choosing $\l_3>0$ and $\l_1,\l_2$ as above
we obtain the strict inequality for the sum of the terms
up to weighted order $8$ for all sufficiently small $(z,u)\ne 0$.
Finally, in the full weighted homogeneous expansion of \eqref{basic},
we will also reach the strict inequality for all sufficiently small $(z,u)\ne 0$
implying $F\in \6J_{p,p'}(D,D')$.
This proves the existence part of Theorem~\ref{nolabel}.
\epf

\appendix
\section{}

\subsection{Attraction property of analytic discs}
The following elementary property has been used in the proof of Lemma~\ref{firstcondbis}
(see \cite{B} for more elaborate refined versions).

\bl\Label{attract} Let $D\subset\C^n$ be a bounded domain and
$p\in\d D$ a boundary point. Suppose that $\1D$ does not contain
nontrivial complex-analytic varieties through $p$. Then, for any
$0<\eta<1$ and any neighborhood $U$ of $p$, there exists another
neighborhood $V$ of $p$ such that, if $f\colon \D\to D$ is a
holomorphic map with $f(0)\in V$, then $f(\eta\D)\subset U$. \el

\bpf By contradiction, suppose that, for some fixed $\eta$ and $U$,
there exists a sequence of holomorphic maps $f_k\colon \D\to D$ with
$f_k(0)\to p$ such that $f_{k}(\eta\D)\not\subset U$. By Montel's
theorem, $\{f_k\}$ can be assumed convergent to a limit map $f\colon
\D\to D$, uniformly on compacta, in particular, on $\eta\1\D$. Since
$f(\D)\subset\1D$ and, by the assumption, $\1D$ does not contain
nontrivial varieties through $p$, we must have $f(z)\equiv p$. The
latter fact implies $f_{k}(\eta\D)\subset U$ contradicting the
choice of the sequence $\{f_k\}$. The proof is complete. \epf

\subsection{Polynomial approximations in real variables}
We begin with a function $f(x)$ in one (real) variable that is
approximated by a polynomial $p(x)$ up to some error term $r(x)$.
We have the following elementary property whose proof is left to
the reader:

\bl\Label{lem-one} Let $p(x)$ be a real polynomial of degree $d$,
$r(x)$ a real function satisfying
$$r(x)=o(|x|^d), \quad x\to 0,$$
and suppose that $p(x)+r(x)\ge 0$ for $x\ge0$ in a neighborhood of
$0$. Then $p(x)\ge 0$ for $x>0$ in a neighborhood of $0$. \el

\br\Label{rem-one} The same statement obviously holds if $d>0$ is
replaced by any real number and $p(x)$ by any finite linear
combination of powers $x^l$ for $l\le d$. \er

We next extend Lemma~\ref{lem-one} to several (real) variables.
For simplicity, we restrict ourselves to the two-dimensional case.
Recall that the {\em Newton polytope} of a polynomial
$p(x_1,x_2)=\sum_{l_1l_2}p_{l_1l_2}x_1^{l_1}x_2^{l_2}$ is the
convex hull of the set of all $(l_1,l_2)$ with $p_{l_1l_2}\ne 0$.
The {\em extended Newton polytope} is the minimal convex set $C$
containing the Newton polytope such that, if $(l_1,l_2)\in C$,
then $(k_1,k_2)\in C$ whenever $k_1\le l_1$ and $k_2\le l_2$. We
have the following extension of Lemma~\ref{lem-one}:

\bl\Label{lem-two} Let $p(x_1,x_2)$ be a real polynomial and for
$j=1,\ldots,s$, let $r_j(x)$ be real functions and $(d_{j1},d_{j2})$ be
pairs of nonnegative integers satisfying
$$r_j(x)=o(|x_1^{d_{j1}}x_2^{d_{j2}}|), \quad x=(x_1,x_2)\to 0, \quad
x_1,x_2\ge 0.$$ Suppose that the convex hull of the set
$\{(d_{j1},d_{j2}):1\le j\le s\}$ does not intersect the interior
of the extended Newton polytope of $p(x)$ and that
$$p(x)+\sum_j r_j(x)\ge 0,$$
for $x_1,x_2\ge 0$ in a neighborhood of $0$. Then $p(x)\ge 0$ for
$x_1,x_2\ge 0$ in a neighborhood of $0$. \el

\bpf It follows from the assumptions that there exists a pair
$(\nu_1,\nu_2)\ne 0$ of nonnegative integers such that, for any
coefficient $p_{l_1l_2}\ne 0$ of $p$ and any $j=1,\ldots,s$, one
has
$$\nu_1 l_1 + \nu_2 l_2 \le \nu_1 d_{j1} + \nu_2 d_{j2}.$$
Then, for any real numbers $\l_1,\l_2>0$, we have $p(\l_1
x^{\nu_1}, \l_2 x^{\nu_2})\ge 0$ for $x>0$ in a neighborhood of
$0$ in view of Remark~\ref{rem-one}. Since $\l_1$, $\l_2$ are
arbitrary, we obtain the conclusion of the lemma. \epf

Consider now the case of variables $X_1\in\R^{n_1}$ and
$X_2\in\R^{n_2}$ and write a polynomial $p(X_1,X_2)$ in the form
$$p(X_1,X_2)=\sum_{l_1l_2}p_{l_1l_2}(X_1,X_2),$$
where $p_{l_1l_2}(X_1,X_2)$ is bihomogeneous in $(X_1,X_2)$ of
bidegree $(l_1,l_2)$. Define the {\em extended bihomogeneous
Newton polytope} of $p$ in $\N^2$ the same way as above. Then, we
obtain the following extension of Lemma~\ref{lem-two}:

\bl\Label{lem-three} Let $p(X_1,X_2)$ be a real polynomial in
$X=(X_1,X_2)\in\R^{n_1}\times\R^{n_2}$ and for $j=1,\ldots,s$, let
$r_j(x)$ be real functions and $(d_{j1},d_{j2})$ pairs of
nonnegative integers satisfying
$$r_j(X)=o(\|X_1\|^{d_{j1}} \|X_2\|^{d_{j2}}), \quad X=(X_1,X_2)\to 0.$$
Suppose that the convex hull of the set $\{(d_{j1},d_{j2}):1\le
j\le s\}$ does not intersect the interior of the extended Newton
polytope of $p(X)$ and that
$$p(X)+\sum_j r_j(X)\ge 0,$$
for $X$ in a neighborhood of $0$. Then $p(X)\ge 0$ for $X$ in a
neighborhood of $0$. \el

The proof can be obtained by restricting $p$ and $r_j$ to the span
of two arbitrary vectors $(v_1,0)$ and $(0,v_2)$ in
$\R^{n_1}\times\R^{n_2}$ and applying Lemma~\ref{lem-two}. In
particular, we have the following ``cancellation rule'' for weighted
homogeneous polynomials:

\bc\Label{weights} Let $\nu_1,\nu_2>0$ be weights assigned to
$X_1,X_2$ and let $p(X_1,X_2)$ be a weighted homogeneous polynomial
of degree $d$ in $(X_1,X_2)\in\R^{n_1}\times\R^{n_2}$, i.e.\
$P(t^{\nu_1} X_1, t^{\nu_2} X_2) = t^d P(X_1,X_2)$ and $r$ be a real
function satisfying
$$r(X_1,X_2)=o((\|X_1\|^{1/\nu_1}+ \|X_2\|^{1/\nu_2})^d), \quad (X_1,X_2)\to
0$$ such that $p(X)+r(X)\ge 0$ for $X=(X_1,X_2)$ in a neighborhood
of $0$. Then $p(X)\ge 0$. Furthermore, if $p_0(X_1,X_2)$ is the
nontrivial bihomogeneous component of $p$ of minimal degree in
$X_1$ (or in $X_2$), then also $p_0(X_1,X_2)\ge 0$. \ec

\subsection{Homogeneous polynomials in complex variables}
By separating homogeneous terms and applying the above statements,
one can reduce general polynomial inequalities to inequalities for
homogeneous terms. We next state some elementary results that can
be useful to separate complex monomials of the form $z^k\bar z^l$.

Let $p(z,\bar z)$ be a homogeneous real-valued polynomial of
degree $d$ with
\begin{equation}\Label{eq-p}
p(z,\bar z)=\sum_{k}p_{k} z^k\bar z^{d-k}\ge 0
\end{equation} for $z\in\C$
in a neighborhood of $0$.

\br \Label{types2} Observe that, if $d$ is odd, then \eqref{eq-p} is
only possible if $p\equiv 0$. If $d$ is even, the situation is more
complicated. Set $d=2s$. By integrating \eqref{eq-p} for
$z=z_0e^{i\theta}$ with $0\le \theta\le 2\pi$, we immediately obtain
that $p_s\ge 0$. \er

In case $p_{s} > 0$ one has, in general, no conclusion about the
other coefficients in \eqref{eq-p}. However, if $p_{s} = 0$, all
other coefficients must vanish:

\bl\Label{types} Let $p(z,\bar z)$ be a homogeneous real-valued
polynomial of degree $2s$ satisfying $\eqref{eq-p}$ for $z\in\C$ in
a neighborhood of $0$. Suppose that $p_{s}=0$. Then $p(z,\bar
z)\equiv 0$. \el

\bpf We assume $p\not\equiv 0$ and prove the statement by induction on the maximal number $k$
with $p_{k}\ne 0$. By the assumption and the reality of $p$, we have $s<k\le 2s$.
Otherwise we have $s<k\le 2s$ and let $\eps$ be any
primitive $4(k-s)$th root of unity. Then, if we multiply $z$ in
\eqref{eq-p} by $\eps$, we obtain a new inequality where the term with $z^k\bar z^{2s-k}$ changes
sign whereas all other terms receive factors different from
$-1$. Hence, by adding the new inequality and the old one, we
eliminate the term with $z^k\bar z^{2s-k}$ and keep all other nonzero
terms with with possibly changed but still nonzero coefficients. By the
induction, the new polynomial must be zero.
This is only possible if $z^k\bar z^{2s-k}$ and its conjugate
are the only nonzero terms of $p(z,\bar z)$. Since $k\ne s$,
we obtain a contradiction with \eqref{eq-p}.
Hence $p(z,\bar z)\equiv0$.
\epf

\end{document}